\documentclass[12pt]{amsart}
\usepackage{graphicx}
\begin{document}
\newtheorem{theorem}{Theorem}[section]
\newtheorem{prop}[theorem]{Proposition}
\newtheorem{lemma}[theorem]{Lemma}
\newtheorem{claim}[theorem]{Claim}
\newtheorem{cor}[theorem]{Corollary}
\newtheorem{defin}[theorem]{Definition}
\newtheorem{example}[theorem]{Example}
\newtheorem{xca}[theorem]{Exercise}
\newcommand{\map}{\mbox{$\rightarrow$}}
\newcommand{\aaa}{\mbox{$\alpha$}}
\newcommand{\Aaa}{\mbox{$\mathcal A$}}
\newcommand{\bbb}{\mbox{$\beta$}}
\newcommand{\ccc}{\mbox{$\mathcal C$}}
\newcommand{\ddd}{\mbox{$\delta$}} 
\newcommand{\Ddd}{\mbox{$\Delta$}}
\newcommand{\Fff}{\mbox{$\mathcal F$}}  
\newcommand{\Ggg}{\mbox{$\Gamma$}}
\newcommand{\ggg}{\mbox{$\gamma$}}
\newcommand{\kkk}{\mbox{$\kappa$}}
\newcommand{\lll}{\mbox{$\lambda$}}
\newcommand{\Lll}{\mbox{$\Lambda$}}
\newcommand{\mlp}{\mbox{$\mu^{+}_{l}$}}
\newcommand{\ml}{\mbox{$\mu_{l}$}}
\newcommand{\mr}{\mbox{$\mu_{r}$}}
\newcommand{\mlpm}{\mbox{$\mu_{l}^{\pm}$}}
\newcommand{\mrpm}{\mbox{$\mu_{r}^{\pm}$}}
\newcommand{\mlm}{\mbox{$\mu_{l}^{-}$}}
\newcommand{\mrp}{\mbox{$\mu_{r}^{+}$}}
\newcommand{\mrm}{\mbox{$\mu_{r}^{-}$}}
\newcommand{\mm}{\mbox{$\mu^-$}}
\newcommand{\mpm}{\mbox{$\mu^{\pm}$}}
\newcommand{\mpp}{\mbox{$\mu^+$}}
\newcommand{\mt}{\mbox{$\mu^{t}$}}
\newcommand{\mb}{\mbox{$\mu_{b}$}}
\newcommand{\mz}{\mbox{$\mu^{\bot}$}}
\newcommand{\mpq}{\mbox{$\mu^{(p,q)}$}}
\newcommand{\omp}{\mbox{$0_{-}^{+}$}}
\newcommand{\oa}{\mbox{$\overline{a}$}}
\newcommand{\ob}{\mbox{$\overline{b}$}}
\newcommand{\opm}{\mbox{$0_{+}^{-}$}}
\newcommand{\opp}{\mbox{$0_{+}^{+}$}}
\newcommand{\Pt}{\mbox{$\tilde{P}$}}
\newcommand{\rrr}{\mbox{$\rho$}} 
\newcommand{\rz}{\mbox{$\rho^{\bot}$}}
\newcommand{\rp}{\mbox{$\rho^+$}}
\newcommand{\rmm}{\mbox{$\rho^-$}}
\newcommand{\rpq}{\mbox{$\rho^{(p,q)}$}}
\newcommand{\sss}{\mbox{$\sigma$}} 
\newcommand{\Sss}{\mbox{$\mathcal S$}} 
\newcommand{\sm}{\mbox{$\sigma^-$}}
\newcommand{\sz}{\mbox{$\sigma^{\bot}$}} 
\newcommand{\spm}{\mbox{$\sigma^{\pm}$}}
\newcommand{\spp}{\mbox{$\sigma^+$}}
\newcommand{\st}{\mbox{$\sigma^{t}$}}
\newcommand{\Ss}{\mbox{$\Sigma$}}
\newcommand{\Th}{\mbox{$\Theta$}} 
\newcommand{\ttt}{\mbox{$\tau$}} 
\newcommand{\bdd}{\mbox{$\partial$}}
\newcommand{\zzz}{\mbox{$\zeta$}}
\newcommand{\qb} {\mbox{$Q_B$}}
\newcommand{\inter}{\mbox{${\rm interior}$}}

\numberwithin{equation}{section}

\title[] {The Goda-Teragaito conjecture: an overview}

\author{Martin Scharlemann}
\address{\hskip-\parindent
        Mathematics Department\\
        University of California\\
        Santa Barbara, CA 93106\\
        USA}
\email{mgscharl@math.ucsb.edu}

\date{\today} 
\thanks{Research supported in part by an NSF grant and RIMS Kyoto}

\begin{abstract} We give an overview of the proof (\cite{Sc}) that the 
only knots that are both tunnel number one and genus one are those 
that are already known: $2$-bridge knots obtained by plumbing together 
two unknotted annuli and the satellite examples classified by 
Eudave-Mu\~noz and by Morimoto-Sakuma.
\end{abstract}
\maketitle

\section{Preliminaries}

\begin{defin}  A graph $\Ggg \subset S^{3}$ is a {\em Heegaard spine} 
if it has a regular neighborhood $H \subset S^{3}$ so that $S^{3} - 
interior(H)$ is a handlebody.
\end{defin}

Note that $\Ggg$ is a Heegaard spine if and only if the decomposition 
$S^{3} = H \cup_{\bdd H} (S^{3} - \inter(H))$ is a Heegaard 
splitting of $S^{3}$. 

\begin{defin} A theta-graph $\theta \subset S^{3}$ is an embedded 
graph consisting of two vertices and three edges, each edge incident to 
both vertices.  A knot $K \subset S^{3}$ has {\em tunnel number one} 
if there is a theta-graph $\theta \subset S^{3}$ so that

\begin{itemize}

\item $\theta$ is a Heegaard spine

\item for some edge $\tau \subset \theta$,  $K = (\theta - 
\inter(\tau))$

\end{itemize}

The edge $\tau$ is called the unknotting tunnel for $K$.  

\end{defin}

\begin{defin}  A knot $K$ has genus one if there is a PL once-punctured 
torus $F \subset S^{3}$ so that $K = \bdd F$.  That is, $K$ has a 
{\em Seifert surface} $F$ of genus one.  
\end{defin}

Both sorts of knots, those of tunnel number one and those of genus 
one, have pleasant and useful properties.  Although each type can be 
quite complicated (as measured, for example, by crossing number), each 
is in some sense the first and easiest collection of knots under one 
natural index of complexity (the tunnel number or the genus).  It's 
therefore of interest to determine which knots are simple in both 
senses.  That is, which knots have both tunnel number one and genus 
one.

The answer, as earlier conjectured by Goda and Teragaito \cite{GT} is 
this:

\begin{theorem} \label{theorem:main} \cite{Sc} Suppose $K \subset 
S^{3}$ has tunnel number one and genus one.  Then either
\begin{enumerate}
\item $K$ is a satellite knot or
\item $K$ is a $2$-bridge knot.
\end{enumerate}
\end{theorem}

\begin{figure} 
\centering
\includegraphics[width=0.15\textwidth]{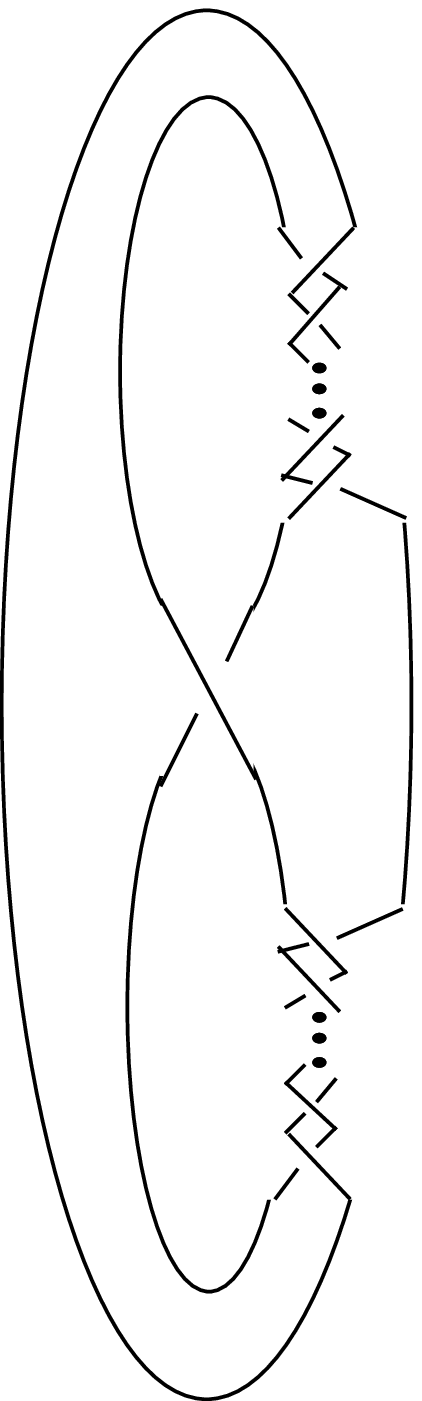}
\caption{} \label{fig:2bridge}
\end{figure}

This theorem is useful because of the historical background: 
$2$-bridge knots all have tunnel number one and those of genus one are 
easily described (cf  \cite[Proposition 12.25]{BZ}) and, indeed, they 
are easily pictured: see Figure \ref{fig:2bridge}.  

On the other hand, Morimoto and Sakuma \cite{MS} and independently 
Eudave-Mu\~noz \cite{EM} classified all satellite knots which have 
tunnel number one.  These knots have a concrete description and can be 
naturally indexed by a $4$-tuple of integers.  In \cite{GT}, Goda and 
Teragaito determined which of these satellite knots have genus one and 
made the conjecture that these knots complete the list of knots that 
have both genus one and tunnel number one.  In other words, they 
conjectured that Theorem \ref{theorem:main} is true.  

The proof that the conjecture is true, as given in \cite{Sc}, is long 
and intricate but much of the energy and length is required by 
arguments that are in some sense technical.  The intention here is to 
give an overview of the proof that focuses on general strategy.  The 
hope is that the reader will understand how the proof follows a rather 
natural course, not one that is as contrived as it might first appear.

\section{The mathematical background} \label{section:math}

Matsuda had verified the Goda-Teragaito conjecture 
(Theorem \ref{theorem:main}) for an important class of knots, and we 
will need his result.  A useful way to state it for our purposes is 
this:

\begin{theorem} \label{theorem:matsuda} \cite{Ma} Suppose that $\tau$ 
is an unknotting tunnel for the genus one knot $K$ and one of the 
cycles in the theta-graph $K \cup \tau$ is the unknot.  Then $K$ is either 
a satellite knot or a $2$-bridge knot.
\end{theorem}

This special case is more general than it might seem.  Suppose, for 
example, that $H$ is a regular neighborhood of a theta-curve Heegaard 
spine $\theta$, so $H$ is a genus two handlebody.  Associated to each 
edge $e$ in $\theta$ is a meridian disk $\mu_{e} \subset H$ that 
intersects $e$ in a single point.  We have this corollary of Matsuda's 
theorem:

\begin{cor} \label{cor:matsuda} Suppose that $K$ is a genus one knot 
lying on $\bdd H$ that intersects $\mu_{e}$ in a single point.  
Suppose further that $\theta - \int(e)$ is the unknot.  Then $K$ is either a 
satellite knot or a $2$-bridge knot.
\end{cor}  

\begin{proof} Apply the ``vacuum cleaner trick'' to the $1$-handle in 
$H$ corresponding to the edge $e$.  That is, think of $K$ as the union of two 
arcs: one arc $\aaa$ runs exactly once over the $1$-handle, and the 
other arc $\bbb$ lies on the boundary of the unknotted solid torus $W 
= H - (1$-$handle)$ and \bbb\ connects the ends of $\aaa$ in $\bdd W$.  Slide 
the ends of the $1$-handle along the arc $\bbb$, vacuuming it up onto 
the $1$-handle until $K$ has been made disjoint from a meridian disk $\mu$ 
of $W$.  At this point, $H$ can be viewed as the regular neighborhood 
of another $\theta$-graph, namely $K \cup \tau$, where $\tau$ is an 
edge intersecting $\mu$ in a single point (i.  e.  $\mu = 
\mu_{\tau}$).  Since $W$ is an unknotted torus, the corresponding cycle in 
$\theta$ is unknotted and Theorem \ref{theorem:matsuda} applies.  
\end{proof}

An unexpected application of Corollary \ref{cor:matsuda} arises from work of 
Eudave-Mu\~noz and Uchida.  Suppose that $H$ is a regular neighborhood 
of a theta curve $\theta$ that is a Heegaard spine.  Let $X = S^{3} - 
\inter(H)$, a genus two handlebody.  Suppose that $F$ is a properly 
imbedded incompressible genus one surface in $X$ with $\bdd F = K = F 
\cap \bdd H$.  

\begin{prop} \label{prop:uchida} Let $A \subset X$ be an 
incompressible annulus obtained from $F$ by $\bdd$-compressing $F$ to 
$\bdd H = \bdd X$.  Suppose there is an edge $e$ of $\theta$ so that 
$|\mu_{e} \cap \bdd A| = 1 = |\mu_{e} \cap K|$.  Then $K$ is either a 
satellite knot or a $2$-bridge knot.
\end{prop}

\begin{proof} $\bdd A$ has two components which we denote $\bdd_{\pm} 
A$ and we take $\bdd_{-} A$ to be the component that intersects 
$\mu_{e}$.  Take two parallel copies of $\mu_{e}$ and band them 
together along the part of $\bdd_- A$ that does not lie between them.  
The result is a disk $E \subset H$ that is disjoint from $\bdd A$ and 
separates $H$, leaving one of $\bdd_{\pm} A$ in the boundary of each 
of the solid tori components of $H - E$.  Label these solid tori 
correspondingly $L_{\pm}$ and denote by $L$ the link whose core 
circles are $L_- \cup L_+$.  Note that $\bdd_- A$ is a longitude of 
$L_-$ and $\bdd_+ A$ is a $(p, q)$ cable of $L_+$, for some $(p, q)$.  

The link $L$ visibly has an incompressible annulus (namely $A$) in its 
complement.  Moreover, $L$ has tunnel number one: attaching to $L$ an edge 
dual to $E$ gives a graph $\Ggg$ whose regular neighborhood is $H$, so 
$\Ggg$ is a Heegaard spine.  Tunnel number one links with essential 
annuli in their complement have been classified by Eudave-Mu\~noz and 
Uchida (cf.  \cite{EU}).  In particular, $L_+$ is the unknot.  But 
$L_{+}$ can also be viewed as the cycle in $\theta$ obtained by 
removing $e$ (equivalently, the core of the torus obtained by removing 
$\mu_{e}$).  Since $|\mu_{e} \cap K| = 1$, the result follows from 
Corollary \ref{cor:matsuda}.
\end{proof}

A special case of this proposition is independently useful:

\begin{cor} \label{cor:uchida} Suppose $\tau$ is an unknotting tunnel 
for $K$ and $\tau$ lies on a genus one Seifert surface $F$ for $K$.  
Then $K$ is either a satellite knot or a $2$-bridge knot.
\end{cor}

\begin{proof} $\tau$ can't be parallel to a subarc of $K$, else $K$ 
would itself be unknotted.  So $\tau$ is an essential arc in $F$, and 
so $F - \eta(\tau)$ is an annulus.  Let $H$ be a regular neighborhood 
of the theta-graph $K \cup \tau$ and $\mu$ be a meridian disk for $K$, 
considered also as a meridian dual to an edge of the graph $K \cup 
\tau$.  Then $F - \inter(H)$ is an incompressible annulus $A$ that 
satisfies the hypotheses of Proposition \ref{prop:uchida}.
\end{proof}

The last proposition and its corollary begin to suggest a strategy for 
trying to prove the conjecture: Let $\theta = K \cup \tau$.  Try to 
arrange that $\tau$ is disjoint from a genus one Seifert surface $F$ 
for $K$, so that we can think of $F$ as lying in the closed complement 
$X$ of $H = \eta(\theta)$ in $S^{3}$.  This makes $F \cap \bdd H$ 
a useful copy of $K$ lying on $\bdd H$.  Try to show that some 
boundary-compressing disk for $F$ in $X$ crosses a meridian of $\tau$ 
exactly once or, alternatively, is disjoint from a meridian 
corresponding to one of the two edges $K - \tau \subset \theta$.  If 
the former happens then, with some work, $\tau$ can be isotoped onto 
$F$ and Corollary \ref{cor:uchida} applies.  If the latter happens 
then we can invoke Proposition \ref{prop:uchida}.

\begin{figure}
\centering
\includegraphics[width=.4\textwidth]{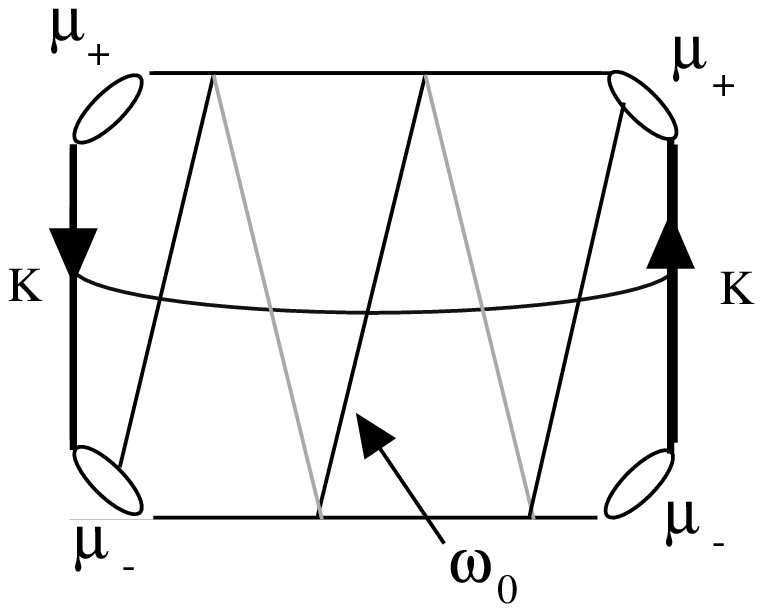}
\caption{} \label{fig:pillowslope}
\end{figure}

Such a strategy requires an understanding of potential 
boundary-compressing disks for $F$ in $X$, once $\tau$ is made 
disjoint from $F$.  A natural source for such disks are the outermost 
disks cut off by $F$ from meridians of the handlebody $X$.  That is, 
if $E$ is a meridian disk of $X$, then an outermost arc of $E \cap F$ 
in $E$ cuts off a disk $E_{0}$.  Moreover, $E_{0}$ lies on one side of 
$F$, so the arc $\omega = E_{0} \cap \bdd H$ has both ends incident to 
the same side of $K$.  Consider the meridians $\mu_{\pm}$ of $H$ 
corresponding to the two edges of $\theta = K \cup \tau$ coming from 
$K$ and the meridian $\mt$ coming from the edge $\tau$.  If $\omega$ 
is disjoint from either of the meridians $\mu_{\pm}$ then we are done 
by Proposition \ref{prop:uchida}.  If $\omega$ intersects both these 
meridians, then some subarc $\omega_{0}$ of $\omega$ cut off by 
$\mu_{\pm}$ is an arc in the $4$-punctured sphere $\bdd H - \mu_{\pm}$ 
running from a copy of $\mu_{-}$ to $\mu_{+}$.  Arcs in a 
$4$-punctured sphere are easy to understand -- each roughly 
corresponds to a rational number given by its slope, viewing the 
$4$-punctured sphere as a pillowcase with holes in the corners.  See 
Figure \ref{fig:pillowslope}
Since $\omega_{0}$ is disjoint from the vertical (i.  e.  slope 
$\infty$) arcs $K \cap \Ss$, it determines a non-zero integral slope 
which (by appropriate choice of slope $0$) we may take to be $1$.  In 
other words, we are done unless the meridian disks of $X$ intersect 
the $4$-punctured sphere $\bdd H - \mu_{\pm}$ in arcs of one 
particular slope.

Of course the genus two handlebody $X$ has many meridians, so it seems 
that it would be difficult to say much about the slopes of arcs 
determined by these meridians.  But it is a remarkable fact (see 
\cite{ST1}) that, if we restrict to simple closed curves in $\bdd H$ 
that bound meridians in both $H$ and $X$, the slope is determined (up 
to only a small ambiguity) by $K$ and $\tau$, as long as $K$ is not 
$2$-bridge.  We are done anyway if $K$ is $2$-bridge, so the upshot is 
that right from the beginning there is only one troublesome case to 
deal with -- when the slope given by a meridian of $X$, a meridian 
whose boundary also bounds a meridian of $E$, is simply $1$.  The 
argument here is more rigorously presented in \cite{ST1} (where what 
we here call the 
``slope'' is there the invariant $\rrr \in \mathbb{Q}/2 \mathbb{Z}$).  

This completes the proof of the Goda-Teragaito conjecture for all but 
the case $\rrr = 1$.  But this final case seems to require a 
substantial broadening of our strategy which we now describe.  
Roughly, we start with $H = \eta(K \cup \tau)$ but find simpler and 
simpler spines for $H$, allowing $K$ to appear more complex on $\bdd 
H$ with respect to these simpler spines.  The hope is to do this in 
such a controlled manner that much of the argument above can still be 
applied and, moreover, the process does not stop until one of the 
cycles in the spine is the unknot.

\section{The psychological background - thin position and its successes} 
\label{section:psych}

The notion of thin position for knots goes back to early work of Gabai 
(\cite{G}).  In outline, the idea is this: Think of $S^{3}$ as swept 
out by an interval's worth of $2$-spheres.  More concretely, choose a 
height function $h: S^{3} \map [-1, 1]$ with exactly two critical 
points: a maximum and a minimum at heights $\pm 1$.  It is possible to 
associate a natural number, called the {\em width}, to any generic 
positioning of a knot $K$ in $S^{3}$.  This can be done so that the 
width is unchanged by isotopies of $K$ that do not create or destroy 
critical points or change the ordering of the heights of the critical 
points.  In fact, width stays unchanged if the height of two adjacent 
maxima or two adjacent minima are switched.  It will go down if a 
maximum is pushed below a minimum or a maximum and a minimum are 
cancelled.  When the width is minimized, the knot is said to be in 
{\em thin position}.  

To illustrate the usefulness of thin position, consider a knot $K 
\subset S^{3}$ in thin position and suppose $F$ is an incompressible 
Seifert surface for $K$.  Suppose some level $2$-sphere $P$ (ie  
$P = h^{-1}(t), t \in (-1, 1)$) is transverse to $F$ and among the 
components of $F - P$ are a pair of disks, say cut off by arc 
components of $P \cap F$, one disk lying above $P$ and one below $P$.  
Then those disks would describe a way to simultaneously push a maximum 
of $K$ below $P$ and a minimum of $K$ above $P$, reducing the width.  
We conclude that no such pair of disks occurs.  

On the other hand, as we imagine level spheres $P_{t} = h^{-1}(t)$ 
sweeping across $F$ from top to bottom, at first there must be a disk 
cut off from $F$ that lies above $P_{t}$ and just before the end of 
the sweep there must be a disk from $F$ that lies below $P_{t}$.  
Since there can't simultaneously be both types of disks, as we have 
just seen, there must be some height at which there are no disks of 
either type, so in fact for $P$ at this height, all components of $P 
\cap F$ are essential in $F$.

If $F$ is of genus one, then essential $1$-manifolds in $F$ are easy 
to describe.  For example, if all components of the $1$-manifold are 
arcs then these arcs lie in at most three different parallelism 
classes.

\bigskip

In \cite{GST} we presented a similar notion of thin position that 
works well for trivalent graphs in $S^{3}$.  That is, just as for 
knots, it is possible to define the width of a generic positioning of 
a trivalent graph $\Ggg$ in $S^{3}$ in such a way that the width is 
unchanged by isotopies of $K$ that do not create or destroy critical 
points or change the ordering of the heights of the critical points, 
where here ``critical point'' is enlarged to include vertices of the 
graph.  Moreover, and this is the delicate point, the width stays 
unchanged if the height of two adjacent maxima or two adjacent minima 
are switched.  Here the delicacy arises because the maxima (and 
minima) can be of two different types: one type is the standard 
maximum of an arc, but the second type is that of a $\lll$-vertex, i.  
e.  a vertex incident to two ends of edges lying below it and one 
above.  The last crucial property of the width is that, just as in the 
case of knots, the width of $\Ggg$ will go down if a maximum is pushed 
below a minimum or a maximum and a minimum are cancelled, or a 
$\lll$-vertex and adjacent standard minimum become a $Y$-vertex (or, 
symmetrically, a $Y$-vertex and adjacent standard maximum become a 
$\lll$-vertex).  To summarize briefly: If $\Ggg$ is in thin position 
then, without affecting width, minima can be pushed past other minima 
and maxima past other maxima but no maximum can be pushed below any 
minimum.

It turns out that any theta-curve Heegaard spine in $S^{3}$ has 
a useful property when it is put in thin position. 

\begin{defin} 
A trivalent graph $\Ggg \subset S^{3}$ is in {\em bridge 
position} with respect to a height function $h: S^{3} \map [-1, +1]$ 
if there is a level sphere $P \subset S^{3}$ so that all maxima of 
$\Ggg$ lie above $P$ and all minima lie below.  Such a sphere $P$ is 
called a {\em dividing sphere}.
\end{defin}

\begin{theorem} \cite{ST2} \label{theorem:levelspine}
Suppose $\theta \subset S^{3}$ is a theta-curve Heegaard spine that is 
in thin position in $S^{3}$.  Then $\theta$ is in bridge position and some 
edge of $\Ggg$ is disjoint from a dividing sphere.
\end{theorem}

With this in mind, we try to extend the application of thin position 
given above to the following: Suppose $\theta$ is a theta-curve 
Heegaard spine with regular neighborhood $H$ and the Seifert surface 
$F$ is properly embedded in $X = S^{3} - \inter(H)$.  Suppose 
$\theta$ is in thin position and $K$ doesn't ``back-track'' along any 
edge in $\theta$.  That is, $K$ doesn't intersect twice in a row the 
meridian corresponding to any one edge of $\theta$.  Then, just as in 
the case for thin position of knots, there is a level sphere $P$ that 
intersects $F$ only in essential curves.  It turns out that the 
following example of this situation (the motivation will come a bit 
later) is crucial:

\begin{defin}  \label{defin:pq}
Suppose $\theta$ is a theta-curve in $S^{3}$ with edges $e^{+}, 
e^{-}, e^{\bot}$.  In $H = \eta(\theta)$, denote the corresponding 
meridians by $\mpp, \mm, \mz$.  Suppose $K \subset \bdd H$ is a 
primitive curve in $\bdd H$ (i.  e.  it intersects some essential disk 
in $H$ in a single point) and $K$ intersects each of the meridians $\mpp, 
\mm, \mz$ always with the same orientation and so that some minimal 
genus Seifert surface $F$ for $K$ intersects $H$ only in $K = \bdd F$.  
Arrange the labelling and orientations of the edges and meridians so 
that, geometrically as well as algebraically,

\begin{itemize}

\item $K \cap \mm = q \geq 1$

\item $K \cap \mpp = p \geq q$

\item $K \cap \mz = p + q$.  

\end{itemize}

Then we say that $K$ (or $(K, F)$) is presented on $\theta$ as a $(p, 
q)$ quasi-cable.

\end{defin}

Our remarks above then show that if $\theta$ is in thin position and 
$(K, F)$ is presented on $\theta$ as a $(p, q)$ quasi-cable then some 
level sphere $P$ intersects $F$ in a $1$-manifold that is essential in 
$F$.  If $F$ is of genus one then the components of $F \cap P$ fall 
into at most three parallelism classes in $F$.  It is not obvious, but 
is not difficult to show, that this greatly constrains the kind of 
$(p, q)$ quasi-cabling that can give rise to genus one surfaces.  Most 
importantly (see \cite[Lemma 2.3]{Sc}) the constraint forces $q = 1$.  
In the context of the strategy outlined above (cf.  Corollary 
\ref{cor:matsuda} and Proposition \ref{prop:uchida}) establishing that 
a meridian dual to an edge intersects $K$ exactly once is a promising 
development indeed.  Indeed, this connection alone suggests that 
thinking about $K$ as a $(p, q)$ quasi-cable could be quite useful.

\bigskip

Next we show why it is not only useful but also perhaps natural to 
think about $K$ as a $(p, q)$ quasi-cable.  For this we go back to the 
argument of Section \ref{section:math} and consider how we might try 
to apply thin position to the only case ($\rrr = 1$) in which the 
argument of Section \ref{section:math} fails.  Recall that $\rrr = 1$ 
means that for $E$ any meridian of $X$, there is an arc $\omega_{0} 
\subset \bdd E$ running between the pair of meridians $\mu_{\pm}$ such 
that $\omega_{0}$ is disjoint from $K$.  Now suppose the theta-graph 
$\theta = K \cup \tau$ is put in thin position and, moreover, among 
all possible ways of sliding $\tau$ on $K$, we've picked the one that 
makes $\theta$ thinnest.  We hope to use the thin position argument 
formerly applied to the Seifert surface $F$ and see what happens if we 
apply it instead to the meridian disk $E$.

A first observation is that, following \cite{GST}, we can assume that 
$K \cup \tau$ is in bridge position and that a dividing sphere is 
disjoint from the tunnel $\tau$.  Indeed it is shown in \cite{GST} 
that, with $K$ in bridge position, the tunnel \ttt\ may be made level 
with its ends at maxima (or minima).  If, when \ttt\ made level, the 
ends of the tunnel are at the same maximum, then the tunnel is 
unknotted and we can appeal to Theorem \ref{theorem:matsuda} to 
complete the argument.  If instead the tunnel \ttt\ has its ends at 
different maxima then the small perturbation that makes $K \cup \tau$ 
generic will leave it in bridge position, with \ttt\ monotonic and 
above the dividing sphere.  

\begin{figure}
\centering
\includegraphics[width=.8\textwidth]{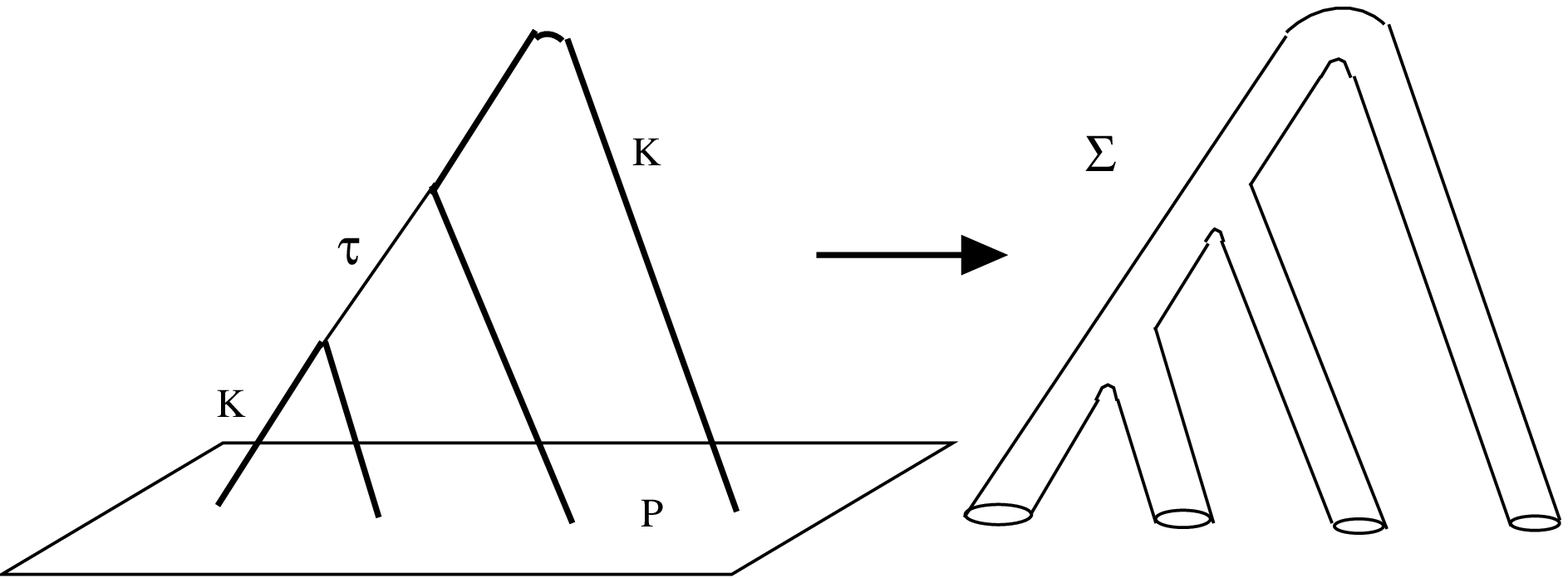}
\caption{} \label{fig:Sigma}
\end{figure} 

\begin{figure}
\centering
\includegraphics[width=.8\textwidth]{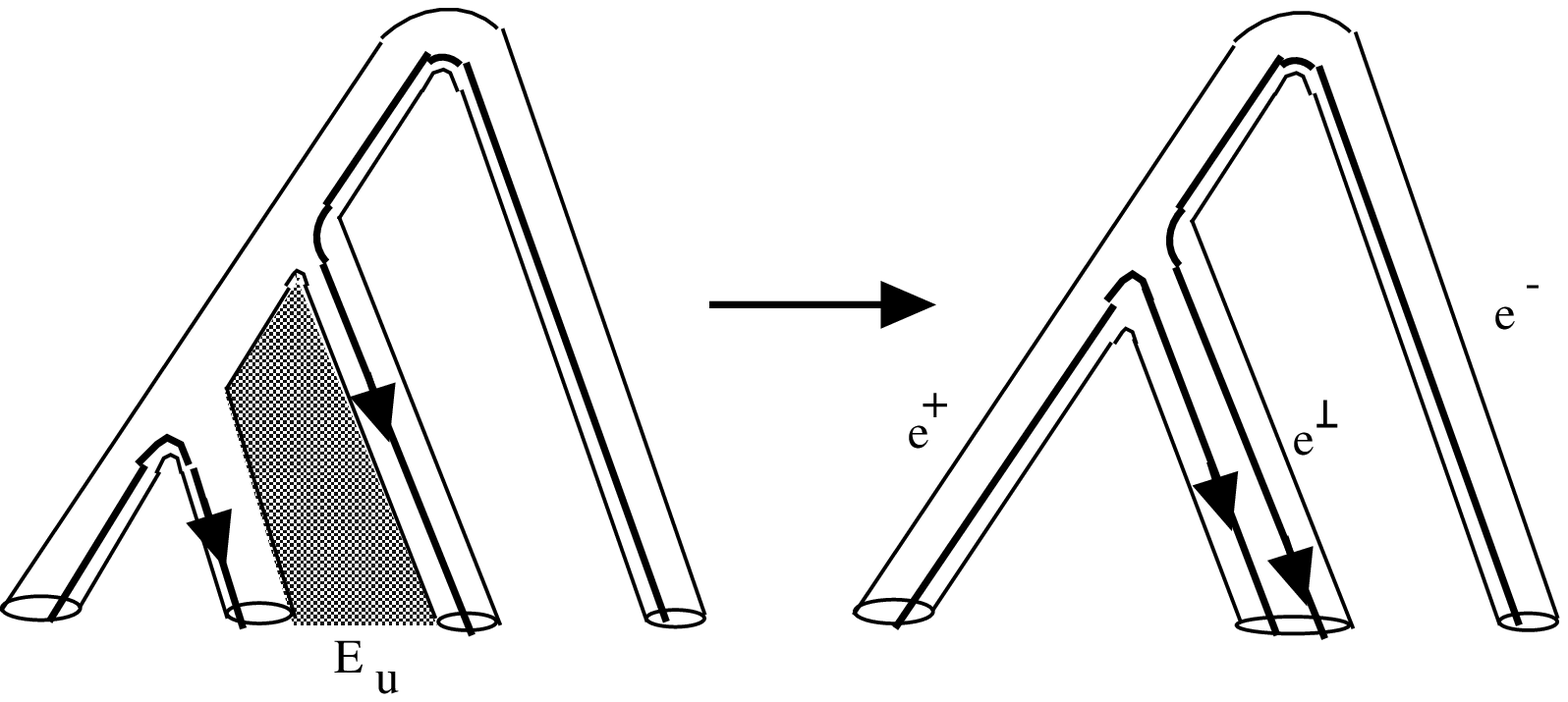}
\caption{} \label{fig:oneone}
\end{figure}

Consider now how level spheres would intersect a meridian disk $E$ of 
$X$ as they sweep across $\theta$.  We saw that as they sweep across 
$F$ there is at least one level at which the intersection with $F$ is 
an essential collection of curves in $F$.  Something must go wrong 
with this argument when we replace $F$ with $E$, for it is impossible 
to have essential arcs in a disk like $E$.  It's easy to see what can 
go wrong: it may be that there is a level sphere $P$ that 
simultaneously cuts off both lower and upper disks from $E$, but these 
disks are incident to $H$ in such a complicated way that they can't be 
used to make $\theta$ thinner.  Such a complicated arc can only arise 
when some component of $\bdd H - P$ is itself complicated, eg the 
$4$-punctured sphere component $\Ss$ of $\bdd H - P$ that appears when 
$P$ is a dividing sphere.  (See Figure \ref{fig:Sigma}.)  In 
particular, we founder when two conditions occur: some subarc 
$\omega_{0}$ of $\bdd E$ is relatively simple in $\Ss$ (so that $\rrr 
= 1$), while another subarc (namely the arc on which an upper disk is 
incident to $\Ss$) is relatively complex.  Since these subarcs are 
disjoint (for $\bdd E$ is embedded) the possibilities are few.  In 
fact, and this is not obvious, the most difficult case to handle is 
one in which the upper disk is so positioned that it defines a way in 
which $\theta$ can only be made thinner by a Whitney move on the edge 
$\tau$.  This Whitney move has the effect of placing $K$ on $\bdd H$ 
as a $(1, 1)$-quasi-cable; see \cite[Lemma 5.5]{Sc} and Figure 
\ref{fig:oneone}.  So we naturally move to thinking about $K$ as a 
quasi-cable on some different spine of $H$.

This Whitney move begins a kind of inductive process.  We try to 
repeat the argument in the more general setting in which $(K, F)$ is 
presented as a quasi-cable on the neighborhood $H$ of a 
theta-curve Heegaard spine.  The knot $K$ is kept thinly presented 
thorughout the process and the strategy is to choose the thinnest 
possible theta-curve for such a presentation.  If such a thinnest 
theta-curve $\theta$ contains an unknotted cycle, then (essentially) 
Corollary \ref{cor:matsuda} can be applied to finish the proof.  If it 
doesn't, then we try to use a meridian $E$ of $X$ to thin $\theta$ 
further and observe that it will do so unless, as previously, $\theta$ 
is in bridge position, some edge of $\theta$ is disjoint from a 
dividing sphere, and some such dividing sphere $P$ cuts off both an upper 
disk $E_{u}$ and a lower disk $E_{l}$ from $E$.  

We then consider again how $\bdd E$ intersects the $4$-punctured 
sphere component $\Ss \subset \bdd H$ of $\bdd H - P$.  If a 
$\bdd$-compressing disk $E_{0}$ for $F$ in $X$ cut off from $E$ by $F$ 
is incident to $\bdd H$ in an arc $\omega$ so short that $\omega$ lies 
entirely in $\Ss$, then the annulus $A$ resulting from the 
$\bdd$-compression naturally satisfies the hypothesis of Proposition 
\ref{prop:uchida}, so we would be done.  The remaining possibility is 
that the arc $\omega \subset \bdd E$ to which $F$ $\bdd$-compresses is 
long.  That is, the subarc $\omega$ of $\bdd E$ is disjoint from $K$, 
has both its ends on the same side of $K$ but passes through $\Ss$ 
perhaps many times.  (We can assume that the ends of $\omega$ both lie 
in $\Ss$.)  The hope would be that this relatively clear picture of 
$\omega$ would be enough to show that the upper and lower disks 
$E_{u}$ and $E_{l}$ could either be used to thin $\theta$ further (a 
contradiction) or at least would describe how to perform a Whitney 
move that would thinly present $K$ as a (possibly more complicated) 
quasi-cable on an even thinner theta-graph.  If this last step works 
(as it did in the simple case $\theta = K \cup \tau)$ then when we 
reach the thinnest possible theta-graph that presents $(K, F)$ as a 
quasi-cable we will be done.  Either a cycle in the theta-graph is 
unknotted and we are done by a version of Corollary \ref{cor:matsuda} 
or the arc $\omega$ lies entirely in $\Ss$ and we are done by 
Proposition \ref{prop:uchida}.

\section{Refinements and combinatorial arguments} 
\label{section:refine}

The last step in the program outlined in Section \ref{section:psych} 
is to exploit upper and lower disks $E_{u}$ and $E_{l}$ cut off from a 
meridian disk $E$ of $X$ by a level sphere $P$.  The hope is that 
$E_{u}$ and $E_{l}$ can be used either to 
thin $\theta$, or at least to perform a useful Whitney move on 
$\theta$.  Just about all we know about the disks $E_{u}$ and $E_{l}$ 
is that they are disjoint from a long arc we've called $\omega$ that is 
itself disjoint from $K$ and has ends incident to $K$ from the same 
side.  We can also assume (from thin position) that $P$ lies below an 
entire edge of $\theta$, so the corresponding component $\Ss$ of $\bdd 
H - P$ is a $4$-punctured sphere.  

\begin{figure}
\centering
\includegraphics[width=.4\textwidth]{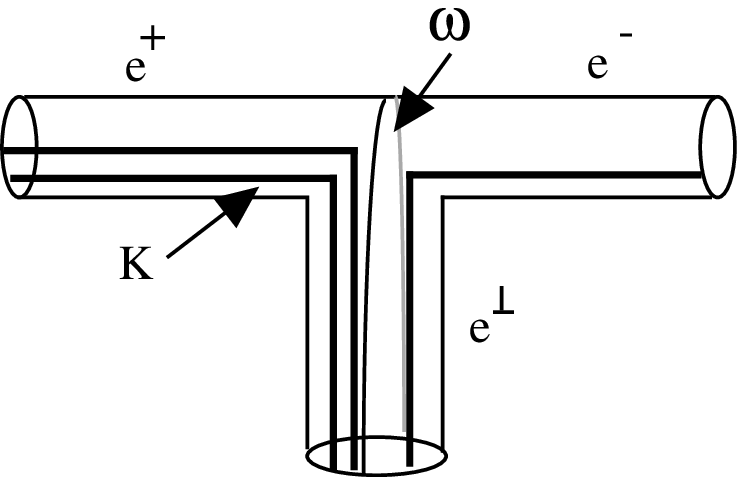}
\caption{} \label{fig:backtrack}
\end{figure}

Things seems rather hopeless, especially if the edge that is disjoint 
from $P$ is $e^{\bot}$: The definition of quasi-cable means that no 
subarc of $K$ goes directly from $e^{+}$ to $e^{-}$.  In particular, 
an arc in $\bdd H$ can ``backtrack'' at an end of $e^{\bot}$ without 
necessarily intersecting $K$.  That is, there is an essential arc in 
$\bdd H - K$ with both ends on $\mz$ and with interior disjoint from 
both of the other meridians $\mpm$.  See Figure \ref{fig:backtrack}.  
This means that, for all we know, the long arc $\omega$ could traverse 
$\Ss$ in a complicated way, passing many times across the meridian 
$\mz$ and still $\omega$ could be disjoint from $K$.  If $\omega$ can 
be this complicated, there is nothing to prevent $\bdd E_{u}$ from 
traversing $\Ss$ in an equally complicated way and so it would be 
useless either for thinning or for describing a simple Whitney move.

Here's an approach to fixing this problem, at least: To ensure that no 
arc of $\bdd E$ backtracks at an end of $e^{\bot}$, somehow guarantee 
that $\bdd E$ {\em does} backtrack at both ends of another edge, one 
of $e^{\pm}$, for then no backtracking arc at an end of $e^{\bot}$ 
could be disjoint from $\bdd E$.  Requiring that $\bdd E$ backtracks 
on one of the other edges $e^{\pm}$ is more reasonable than it might 
appear: Recall that $E$ could have been chosen so that its boundary also 
bounds a disk $D \subset H$.  If we consider how $D$ intersects the 
meridians $\mpm, \mz$ of $H$, an outermost arc of intersection on $D$ 
will precisely cut off a backtracking arc from $\bdd D = \bdd E$.  So 
we know that there is a back-tracking arc (called a {\em wave}) on 
exactly one edge.  So the refinement needed to avoid backtracking at 
an end of $e^{\bot}$ is to insist that we will only consider 
theta-graphs which present $(K, F)$ as a quasi-cable and which also 
satisfies the {\em wave condition}; that is, we require that there is 
a wave of $\bdd E$ based at one of the meridians $e^{\pm}$.  Of 
course, this restriction means that the inductive step (i.  e.  the 
Whitney move) must be shown to preserve the wave condition at the next 
stage, including the first stage in which $K \cup \tau$ becomes a 
theta graph presenting $(K, F)$ as a $(1, 1)$ cable.  Checking this is 
a bit technical and we only note here that it is shown in \cite{Sc} 
that it can be done.

Introducing the wave condition gives unexpectedly useful information 
about the long arc $\omega$.  If, for example, a wave is based at the 
meridian $\mpp$ of $e^{+}$ then $\omega$ can't pass directly from an end 
of $e^{-}$ to an end of $e^{\bot}$ since a wave is in the way, nor can 
it backtrack on any edge, for even in backtracking at an end of 
$e^{+}$ it would cross $K$.  It follows, for example, that $\omega$ 
always crosses each meridian $\mpm, \mz$ in the same direction.  Even 
more (somewhat technical) information is available about $\omega$ and 
hence about the boundary of the annulus $A$ that is obtained from $F$ 
by $\bdd$-compressing to $\omega$, but what we've described is enough to 
give an outline of the rest of the argument.  

The first observation is that, if $\aaa_{u} = E_{u} \cap \bdd H$ were 
actually parallel to a segment of $\omega \cap \Ss$ we would be nearly 
finished: Since $\omega$ doesn't backtrack along any edge, there are 
only a few paths it can take through $\Ss$; if $\aaa_{u}$ also took 
one of these paths, then it's relatively easy to show that $E_{u}$ and 
$E_{l}$ describe a way either to thin $\theta$ or to perform an 
appropriate Whitney move (depending on the path).  It turns out that, 
when the edge of $\theta$ that is disjoint from $P$ is the edge that 
contains the wave, just knowing that the wave is disjoint from $\omega$ 
suffices to show there is an appropriate Whitney move.  So, for a good 
representative of the hard cases that remain, we assume 
that it is the edge $e^{\bot}$ that is disjoint from the dividing 
sphere $P$.

\bigskip

Now a brief digression.  We've noted above that if $\aaa_{u}$ were 
parallel to a segment of $\omega \cap \Ss$ we would be finished.  This 
suggests (and one can prove) that there is a dividing sphere $P$ that 
intersects the annulus $A$ only in essential spanning arcs.  Roughly, 
the argument is the standard thin position argument: if, as $P_{t}$ 
sweeps across $\theta$, there were always some inessential arc of 
intersection then at some level there would be a dividing sphere that 
cut off both an upper and a lower disk from $A$.  But $\bdd A$ does 
follow (roughly) $\omega$, since $\bdd A$ is obtained by banding $K$ 
to itself along $\omega$.  Now these two disks could be used, either 
to thin $\theta$ or to perform an appropriate Whitney move, changing 
$\theta$ to a thinner graph which still presents $(K, F)$ as a 
quasi-cable.

\bigskip

In this last hard case it appears to be difficult to complete the 
argument as planned, using $E_{u}$ and $E_{l}$, simply because the 
information that we have (i.  e.  that $\aaa_{u}$ is disjoint from 
$\omega$) doesn't seem quite enough to describe a good Whitney move.  
Indeed, this information isn't even enough to show that $E_{u}$ is 
disjoint from $K$ so it could hardly be used to describe a Whitney 
move on $\theta$ that wouldn't disrupt the presentation of $K$.  But 
in the digression and earlier we have seen that in this hard case the 
annulus $A$ to which $F$ $\bdd$-compresses has some nice properties: 
Its boundary is ``regular'' (i.  e.  there is no back-tracking) and 
it's relatively easy to describe the type of path its boundary 
follows, since each boundary component of $A$ consists of a copy of 
$\omega$ and part of $K$.  Moreover, $A \cap P$ consists exclusively 
of parallel spanning arcs in $A$.  This suggests that a combinatorial 
approach might be useful, as at the beginning when we used a 
description of the essential collection of curves $F \cap P \subset F$ 
to establish that $q = 1$.

In fact the same sort of combinatorial argument does allow us to 
simplify even further the possible paths that $\bdd A$ can follow in 
$\bdd H$, but again the argument seems to stop just short of 
completing the proof.  The exceptional case is quite specific.  Each 
component of $\bdd A$ can be described as a path in $\bdd H$ simply by 
writing down a word in letters $a, b, \oa, \ob$ that describe the 
order in which the curve intersects $\mpp$ and $\mm$.  If it 
intersects $\mm$ with the same orientation that $K$ does, write down 
$a$; it if intersects $\mpp$ with the opposite orientation that $K$ 
does, write down $\ob$.  Then it turns out that a combinatorial 
argument, using the parallel arcs of intersection of $A \cap P$ in 
$A$, eliminates all possibilities for $\omega$ except a word $w$ that 
is positive in the letters $a$ and $\ob$.  It also eliminates all 
possibilities for $\bdd A$ except the words

\bigskip

\begin{itemize}
\item $\bdd_- A \leftrightarrow wa$
\item $\bdd_+ A \leftrightarrow w\ob^{p}$.
\end{itemize}

\bigskip

The remarkable regularity of this structure opens another possible way 
of eliminating this last remaining possibility: use a combinatorial 
argument in the graph formed in $P$ by the arcs $A \cap P$.  That is, 
consider the graph $\Lll$ formed in $P$ by letting intersection points 
of $e^{\pm}$ with $P$ be the vertices and the arc components $A \cap 
P$ be the edges.  (The idea of using such graphs seems to go back at 
least to \cite{GL}.)  Remarkable order can be discerned.  There are no 
inessential loops.  If we orient the edges in $P$ so that the 
corresponding orientation in $A$ is from $\bdd_{+} A$ to $\bdd_{-} A$, 
there cannot be two faces of $P - \Lll$ (from the same component of 
$\Lll$) that have their boundaries consistently oriented in opposite 
directions, one clockwise and the other counterclockwise.  In any case 
(suppressing technicalities) the argument now proceeds by showing that 
such opposite pairs of faces must always arise, first by showing that 
they arise when $w$ is long, then by successively eliminating all 
shorter possibilities for $w$.

\section{Epilogue - Heegaard spines that are eyeglass graphs}

The argument above is only a sketch, but a particularly important 
possibility has been glossed over to smooth the flow of the argument.  
This final section gives a brief description of how the case arises 
and how it is handled.

\begin{figure}
\centering
\includegraphics[width=.6\textwidth]{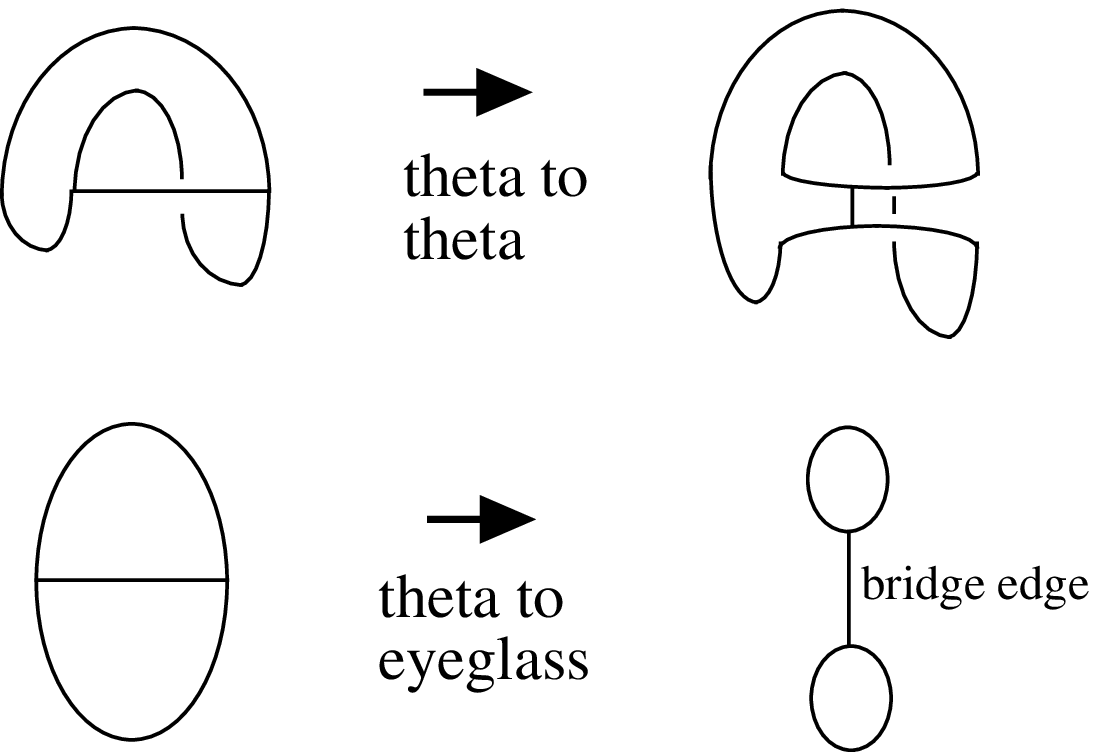}
\caption{} \label{fig:thetaeye}
\end{figure}

As the argument has been described, thin position puts $\theta$ in 
bridge position with an edge disjoint from a dividing sphere.  Then a 
surface in $X$ (e.  g.  the annulus $A$) is used to find disjoint 
upper and lower disks.  These disks are used to describe a Whitney 
move on $\theta$, changing it to an even thinner theta-graph which 
still thinly presents $(K, F)$ as a quasi-cable.  What we have 
suppressed is the possibility that the Whitney move described by the 
disks changes $\theta$ not to another theta-graph but instead to an 
``eyeglass'' graph -- that is, the graph formed by connecting two 
loops on different vertices with an edge, called the bridge edge. See 
Figure \ref{fig:thetaeye}.  In 
this epilogue we briefly describe what to do in this case.

The first important point is that there is a good theory of thin 
eyeglass Heegaard spines in $S^{3}$, as there is for theta-graphs, cf 
\cite{ST2}.  In particular, when an eyeglass Heegaard spine $\bowtie$ 
is put in thin position, either it's in bridge position with the 
bridge edge disjoint from a dividing sphere, or at least one of the 
two loops in $\bowtie$ is unknotted.  For our purposes, some further argument is needed 
to show that the moves used to thin $\bowtie$ would not destroy the 
thinness of a knot $K$ sitting appropriately on the boundary of its 
regular neighborhood $H$.  The way in which $K$ sits on $\bdd H = 
\bdd\eta(\bowtie)$ is easy to describe (it's called a {\em 
$p$-eyeglass presentation} in \cite[Definition 4.1]{Sc}) and so it is 
straightforward to verify that the thinning moves on $\bowtie$ do not 
disrupt the thinness of $K$.  If, once $\bowtie$ is thin, one of the 
loops is the unknot, the argument concludes essentially via Theorem 
\ref{theorem:matsuda}.  If instead all that we know is that the bridge 
edge is disjoint from a dividing sphere, it's still easy to verify 
also that there is an associated embedding of $\theta$ that is no 
thicker than $\bowtie$.

The upshot is this: If a $p$-eyeglass presentation is needed because 
the wrong sort of Whitney move is given by the upper and lower disks, 
then thin position can be applied to the eyeglass curve $\bowtie$ as 
well, in a way that leads to a contradiction.  That is, if a Whitney 
move converts a theta-graph $\theta$ into a $p$-eyeglass presentation 
on $\bowtie$, then $\bowtie$, as imbedded, is thinner than $\theta$ 
was.  Now isotope $\bowtie$ so it is in thin position.  This can be 
done without affecting the fact that $K$ is in thin position and, 
afterwards, $\theta$ itself could be recovered from $\bowtie$ in a way 
that maintains the same width as $\bowtie$, and so actually leaves 
$\theta$ thinner than when it begain.  This is a contradiction, since 
$\theta$ was originally supposed to be in thin position.  So the role 
of the $p$-eyeglass is merely ``catalytic'' in the sense that the 
eyeglass curve $\bowtie$ appears briefly, and only so that thin 
position applied to $\bowtie$ carries the original imbedding of 
$\theta$ to an even thinner imbedding of $\theta$.  But it's 
appearance is consistent with and indeed reinforces the central theme 
of the argument: thin position for Heegaard spines is a useful tool to 
understand the geometry and topology of knots that lie on their 
corresponding Heegaard surfaces.

\end{document}